\documentclass{article}

\usepackage{amsmath,amsthm,amssymb}

\usepackage{pgfplots}
\pgfplotsset{compat=1.18}
\usepackage{subcaption}
\usepackage{color}
\usepackage{soul}

\newtheorem{theorem}{Theorem}
\newtheorem{lemma}{Lemma}
\newtheorem{corollary}{Corollary}
\newtheorem{proposition}{Proposition}
\newtheorem{conjecture}{Conjecture}

\makeatletter
\renewcommand{\@fnsymbol}[1]{\@alph{#1}}
\makeatother

\begin{document}

\title{Real-Analytic Capacity of Circular Perturbations}
\author{Stephen Deterding \thanks{Marshall University} \and Dante Mancino \thanks{Bucknell University} \and Aiden Taylor \footnotemark[1] \and Lincoln Williston \thanks{Calvin University}}
\maketitle

\begin{abstract}
    The goal of this paper is to prove that the analytic capacity of a family of circular perturbations with boundary $e^{it} + ae^{ikt}: t \in[0,2\pi),\,a\in(-1/|k|,1/|k|),\,k\in \mathbb{Z}$ is real-analytic with respect to $a$ near $a=0$. For $k<0$ we show $\gamma(e^{it} + ae^{ikt})$ is constantly equal to 1, and the same is trivially true at $k=0$. At $k=1$ the capacity is equal to $1+a$. At $k=2$, we write an explicit formula for the  capacity as an infinite product in terms of a (using the Ramanujan Psi function). Furthermore, when $k=2$ we explicitly solve for the analytic capacity at $a=\pm\frac{1}{2}$, which turns out to be $\pi^2/8$.
\end{abstract}

\begin{figure}[htbp]
\centering

\begin{subfigure}[b]{0.31\textwidth}
\centering
\begin{tikzpicture}
\begin{axis}[
    width=\textwidth,
    axis equal image,
    xmin=-1.6, xmax=1.6,
    ymin=-1.6, ymax=1.6,
    axis lines=middle,
    xtick={-1,0,1},
    ytick={-1,0,1},
    ticklabel style={font=\tiny},
    axis line style={gray!60},
    tick style={gray!60},
]
\addplot[domain=0:360, samples=200, dashed, gray, thin]
    ({cos(x)}, {sin(x)});
\addplot[domain=0:360, samples=400, thick, red!65!black]
    ({cos(x)+0.5*cos(2*x)}, {sin(x)+0.5*sin(2*x)});
\end{axis}
\end{tikzpicture}
\caption{$k=2,\ a=0.5$ (cusp)}
\end{subfigure}
\hfill
\begin{subfigure}[b]{0.31\textwidth}
\centering
\begin{tikzpicture}
\begin{axis}[
    width=\textwidth,
    axis equal image,
    xmin=-1.6, xmax=1.6,
    ymin=-1.6, ymax=1.6,
    axis lines=middle,
    xtick={-1,0,1},
    ytick={-1,0,1},
    ticklabel style={font=\tiny},
    axis line style={gray!60},
    tick style={gray!60},
]
\addplot[domain=0:360, samples=200, dashed, gray, thin]
    ({cos(x)}, {sin(x)});
\addplot[domain=0:360, samples=400, thick, red!65!black]
    ({cos(x)+0.2*cos(3*x)}, {sin(x)+0.2*sin(3*x)});
\end{axis}
\end{tikzpicture}
\caption{$k=3,\ a=0.2$}
\end{subfigure}
\hfill
\begin{subfigure}[b]{0.31\textwidth}
\centering
\begin{tikzpicture}
\begin{axis}[
    width=\textwidth,
    axis equal image,
    xmin=-1.6, xmax=1.6,
    ymin=-1.6, ymax=1.6,
    axis lines=middle,
    xtick={-1,0,1},
    ytick={-1,0,1},
    ticklabel style={font=\tiny},
    axis line style={gray!60},
    tick style={gray!60},
]
\addplot[domain=0:360, samples=200, dashed, gray, thin]
    ({cos(x)}, {sin(x)});
\addplot[domain=0:360, samples=400, thick, red!65!black]
    ({cos(x)+0.18*cos(5*x)}, {sin(x)+0.18*sin(5*x)});
\end{axis}
\end{tikzpicture}
\caption{$k=5,\ a=0.18$}
\end{subfigure}

\caption{The curve $K_a$ for the unit circle ($\alpha=1$, dashed) under increasing
parameter $k$. Each $a$ is chosen relative to the sharp threshold $\alpha/k$ of
Proposition~\ref{prop:sharp}: exactly at the threshold for $k=2$, where the curve develops a
cusp at $t=\pi$, and comfortably below it for $k=3,5$.}
\label{fig:perturbations}
\end{figure}
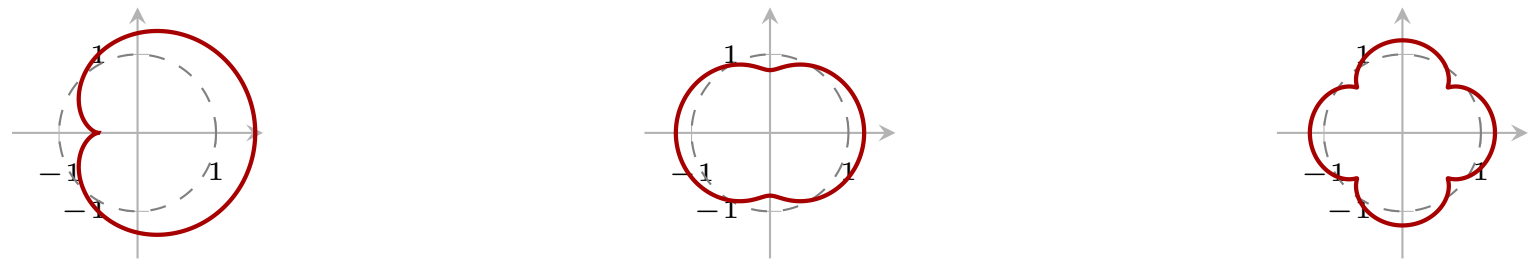

\section{Introduction}

The Painlevé problem is the desire to discover what characterizes removable sets. If for every bounded function $f$ analytic on an open set $\Omega$ with another open set $E\in \Omega$, $f$ is also analytic on the open set $\Omega\setminus E$, then $E$ is a removable set. Ahlfors resolved this decades later by introducing analytic capacity, and he showed that a set is removable if and only if its analytic capacity is equal to 0. The analytic capacity of a set E is defined as the supremum of $|f'(\infty)|$ taken over admissible functions of $E$. A function is considered admissible for a set E if it is both bounded ($|f(z)|\le 1$) and analytic over $\hat{\mathbb{C}}\setminus E$, the complement of E in the extended complex plane, with the additional criterion that $f(\infty)=0$. Analytic capacity has uses outside the Painlevé problem, especially in rational approximation theory. Vitushkin utilizes analytic capacity as part of a characterization of when every continuous function on a set can be uniformly approximated by rational functions with poles outside the set \cite{vitushkin1967}. It is also used to characterize peak points for $R(X)$, the uniform closure of rational functions with poles off $X$ \cite{melnikov1966}.

Analytic capacity, despite its abstract formalization, has some nice, intuitive properties. It is rotationally invariant, scales linearly, and has closed-form solutions for geometrical sets like ellipses and polygons. Outside of these shapes, very few closed-form solutions are known. The analytic capacity of any set and the capacity of its boundary are equivalent. The analytic capacity and logarithmic capacity of Jordan curves are identical, a fact needed for much of the proof of analyticity. This enables us to write 'capacity' as shorthand for either analytic or logarithmic capacity in case of the curves discussed below. 

We now turn our attention to circular perturbations such as those in Figure~\ref{fig:perturbations}. By calculating the capacity of these perturbations, it may be possible to solve for the capacities of other shapes that are heretofore unknown.
In section 2, we proceed with proving the analyticity near $a=0$ of the circular perturbations parametrized by $e^{it}+ae^{ikt}$, where $t$ ranges from $0$ to $2\pi$, $a$ ranges from $0$ to $1/|k|$, and $k$ is an integer. In section 3, we begin by calculating the coefficients for the simpler cases of $k\le 1$. The case $k=2$ admits special attention for multiple reasons. First, we derive a closed-form solution for $\gamma(K_a)$ as an infinite product in terms of $a$. Next, we manipulate the closed-form to solve for the capacity at the boundary $a=\pm1/2$, which turns out to be $\pi^2/8$. Finally, we prove all of the $k=2$ taylor coefficients of $\gamma(K_a)$ are integers. We conclude by conjecturing this pattern holds for all $k\ge2$.

\section{Analyticity of the Perturbed Circle Family}

Fix $\alpha>0$, $k\in\mathbb Z$, $m:=k-1$, and for real $t$ and $|a|<\alpha$ write
\[
z_a(t)=\alpha e^{it}+a e^{ikt}=e^{it}h_a(t),
\qquad h_a(t):=\alpha+a e^{imt},
\]
and, for real $a$ with $|a|$ small, let $K_a:=z_a([0,2\pi])$,
let $P_a$ be the closed Jordan region it bounds, and
$\Omega_a:=\widehat{\mathbb C}\setminus P_a$. We write
$\operatorname{cap}(\cdot)$ for logarithmic capacity, reserving $\gamma_j$ for the Taylor
coefficients of $a\mapsto\operatorname{cap}(K_a)$. On a continuum this
agrees with analytic capacity\cite{ahlfors1946,garnett1972}, and since $K_a$ and $P_a$ share
the unbounded complementary component $\Omega_a$, the capacity depends only on
$\Omega_a$. Our goal is to prove that $a\mapsto\operatorname{cap}(K_a)$ is
real-analytic near $0$.

The argument runs in four movements. Step~A puts the curve in polar coordinates and produces a
single real-analytic function $\Lambda_a$ (the polar log-radius function),
holomorphic in $a$. Step~B shows that composing $\Lambda_a$ with an unknown
boundary reparametrization is a well-behaved, jointly holomorphic operation on a Banach algebra
of analytic functions. Step~C solves for that reparametrization as a contraction fixed point,
builds the exterior conformal map explicitly from it, and confirms via the argument principle
that the construction is the genuine exterior map. Step~D reads off analyticity of the capacity.

\subsection{Notation and constants}
\label{subsec:notation}

Throughout, $S_\eta:=\{t\in\mathbb C:|\operatorname{Im}t|<\eta\}$ denotes the horizontal strip
of half-width $\eta$, and $\mathbb T=\mathbb R/2\pi\mathbb Z$. We collect the standing objects
and constants for reference.

\emph{Polar data (holomorphic continuation).} Since $\operatorname{Im}$ and $|\cdot|$ are not
holomorphic, we continue the polar decomposition $z_a=\rho_a e^{i\theta_a}$
off the real axis by introducing the reflected factor
\[
h^{\#}_a(t):=\alpha+a e^{-imt}
\quad(\,=\overline{h_a}\text{ on }\mathbb R\text{ for real }a\,),
\qquad
\rho_a:=\bigl(h_a h^{\#}_a\bigr)^{1/2},
\]
\[
\theta_a(t):=t+\frac{1}{2i}\bigl(\operatorname{Log}h_a(t)-\operatorname{Log}h^{\#}_a(t)\bigr)
\tag{1}
\]
using the principal branches of the logarithm. On $\mathbb R$, for real $a$, these agree with $|z_a|$
and $\arg z_a$; by the identity theorem they are the unique holomorphic continuations,
$2\pi$-periodic in the sense $\rho_a(t+2\pi)=\rho_a(t)$,
$\theta_a(t+2\pi)=\theta_a(t)+2\pi$. We write
$\Lambda_a:=\log\rho_a\circ\theta_a^{-1}$ (the log-radius as a
function of angle) and $\Xi_a:=\Lambda_a-\log\alpha$.

\emph{The Banach algebra $\mathcal A_\mu$.} For $\mu>0$, let $\mathcal A_\mu$ be the space of
$2\pi$-periodic $f=\sum_n\hat f(n)e^{in\sigma}$ with
$\|f\|_\mu:=\sum_n|\hat f(n)|e^{\mu|n|}<\infty$. The map $\Phi:f\mapsto(\hat f(n)e^{\mu|n|})_{n\in\mathbb Z}$ is a linear
bijection from $\mathcal A_\mu$ onto $\ell^1(\mathbb Z)$, and by definition of $\|\cdot\|_\mu$
it is isometric: $\|f\|_\mu=\|\Phi(f)\|_{\ell^1}$. Since $\ell^1(\mathbb Z)$ is complete and
completeness transfers along isometric isomorphisms, $\mathcal A_\mu$ is complete.
 Using $e^{\mu|n|}\le e^{\mu|k|}e^{\mu|n-k|}$, which follows from the triangle
inequality $|n|\le|k|+|n-k|$, we get
\[
\|fg\|_\mu\le\sum_n\sum_k|\hat f(k)||\hat g(n-k)|e^{\mu|k|}e^{\mu|n-k|}
=\Big(\sum_k|\hat f(k)|e^{\mu|k|}\Big)\Big(\sum_m|\hat g(m)|e^{\mu|m|}\Big)=\|f\|_\mu\|g\|_\mu,
\]
after reindexing $m=n-k$ and applying Fubini's theorem followed by Tonelli's theorem for the (absolutely convergent) double
sum. Clearly $\|e^{in\sigma}\|_\mu=e^{\mu|n|}$, and since $e^{\mu|n|}\ge1$,
\[
\|f\|_\infty\le\sum_n|\hat f(n)|\le
\sum_n|\hat f(n)|e^{\mu|n|}=\|f\|_\mu. \tag{2}
\] 

For $\sigma+i\tau\in S_\mu=\{|\mathrm{Im}\,\zeta|<\mu\}$, the
series $f(\sigma+i\tau)=\sum_n\hat f(n)e^{in\sigma}e^{-n\tau}$ satisfies
$|\hat f(n)e^{in\sigma}e^{-n\tau}|\le|\hat f(n)|e^{\mu|n|}$ for $|\tau|\le\mu$, so by the
Weierstrass $M$-test it converges absolutely and uniformly on every closed sub-strip
$|\mathrm{Im}\,\zeta|\le\mu-\epsilon$; being a locally uniform limit of entire functions
$\hat f(n)e^{in\zeta}$, it is holomorphic on $S_\mu$.

Termwise differentiation (justified by the same uniform convergence)
gives $f'(\sigma)=\sum_n in\hat f(n)e^{in\sigma}$, so
\[
\|f'\|_\infty\le\sum_n|n|\,|\hat f(n)|=\sum_n\big(|n|e^{-\mu|n|}\big)|\hat f(n)|e^{\mu|n|}
\le\Big(\sup_{x\ge0}xe^{-\mu x}\Big)\|f\|_\mu=\frac1{e\mu}\|f\|_\mu, \tag{3}
\]
using $\sup_{x\ge0}xe^{-\mu x}=1/(e\mu)$ (attained at $x=1/\mu$).

\emph{The Hilbert transform.} $\widehat{\mathcal Hf}(n):=-i\operatorname{sgn}(n)\hat f(n)$; thus
$\mathcal H$ annihilates constants, and its multiplier has modulus $\le1$ with weights
depending only on $|n|$, so $\|\mathcal H\|_{\mathcal A_\mu}\le1$. If $f$ is real-valued then
$\hat f(-n)=\overline{\hat f(n)}$, whence $\widehat{\mathcal Hf}(-n)=\overline{\widehat{\mathcal Hf}(n)}$
and $\mathcal Hf$ is real-valued.

\emph{Constants.} With $m=k-1$,
\[
\delta:=\frac{\log2}{m},\quad
\Sigma_0:=\sum_{n\in\mathbb Z}e^{-\frac3{16}\delta|n|}<\infty,  \quad
a_*:=\min\Bigl\{\tfrac{\alpha}{8(k+1)},\ \tfrac{\alpha\delta}{32}, \tfrac{\alpha\delta}{640\,\Sigma_0}\Bigr\}
\]
and $D_*:=\{|a|<a_*\}$. Note
$a_*<\alpha/k$, and $a_*$ is chosen so that
$\frac{10|a|}{\alpha}\Sigma_0\le\delta/64$ on $D_*$. On $S_\delta$ we have
$|a e^{\pm imt}|\le2|a|$, so for $a\in D_*$
\[
|h_a-\alpha|\le2|a|\le\tfrac{\alpha}{12}
\tag{4}
\]
(and likewise for $h^{\#}_a$); both lie in $D(\alpha,\alpha/12)$, off the branch cut,
so (1) is holomorphic on $S_\delta\times D_*$.

\subsection{Preliminary lemmas}

We record two Banach-space lemmas used repeatedly; standard facts (Cauchy estimates, existence of a holomorphic logarithm, and the fact that a bounded linear map
sends holomorphic Banach-valued maps to holomorphic ones) are used without further comment\cite{mujica1986,ahlfors1979}.

\begin{lemma}[Fourier criterion for holomorphy]
\label{lem:fourier-hol}
Let $0<\mu<\nu$ and $D\subset\mathbb C^d$ open, and for $p\in D$ let $f_p$ have Fourier
coefficients $c_n(p)$. If each $c_n$ is holomorphic on $D$ and $|c_n(p)|\le We^{-\nu|n|}$ for
all $n,p$, then $p\mapsto f_p$ is holomorphic from $D$ into $\mathcal A_\mu$, with
$\sup_D\|f_p\|_\mu\le W\Sigma(\nu-\mu)$, where $\Sigma(\beta):=\sum_{n\in\mathbb Z}e^{-\beta|n|}$.
\end{lemma}

\begin{proof}
The partial sums $S_N(p)=\sum_{|n|\le N}c_n(p)e^{in\sigma}$ are holomorphic (finite sums of
scalar holomorphic functions times fixed vectors), and
$\|f_p-S_N(p)\|_\mu\le W\sum_{|n|>N}e^{-(\nu-\mu)|n|}\to0$ uniformly on $D$; the case $N=0$
gives the norm bound. A locally uniform limit of holomorphic $\mathcal A_\mu$-valued maps is
holomorphic\cite{mujica1986}.
\end{proof}

\begin{lemma}[Holomorphic dependence of contraction fixed points]
\label{lem:fixed-point}
Let $X$ be a complex Banach space, $x_0\in X$, $\rho>0$, and
$\overline B_\rho:=\{x:\|x-x_0\|\le\rho\}$. Let $D\subset\mathbb C^d$ be open and $U\supset\overline B_\rho$
open in $X$, and let $T:U\times D\to X$ be jointly holomorphic with
$T(\overline B_\rho\times D)\subset\overline B_\rho$ and
$\|T(x,p)-T(y,p)\|\le\kappa\|x-y\|$ for all $x,y\in\overline B_\rho$, $p\in D$, with $\kappa<1$.
Then $T(\cdot,p)$ has a unique fixed point $x(p)\in\overline B_\rho$, the map $p\mapsto x(p)$ is
holomorphic on $D$, and
\[
\|x(p)-x_0\|\le\frac{\|T(x_0,p)-x_0\|}{1-\kappa}.
\]
If moreover $X_{\mathbb R}\subset X$ is a closed real subspace with $x_0\in X_{\mathbb R}$,
$D_{\mathbb R}:=D\cap\mathbb R^d$, and
$T\bigl((\overline B_\rho\cap X_{\mathbb R})\times D_{\mathbb R}\bigr)\subset X_{\mathbb R}$, then
$x(p)\in X_{\mathbb R}$ for $p\in D_{\mathbb R}$.
\end{lemma}

\begin{proof}
Existence and uniqueness are the Banach fixed point theorem\cite{rudin1991} on the complete metric space
$\overline B_\rho$. The Picard iterates $x^{(0)}\equiv x_0$, $x^{(j+1)}(p):=T(x^{(j)}(p),p)$ are
holomorphic (compositions of $T$ with holomorphic maps into $U$) and satisfy
$\|x^{(j)}(p)-x(p)\|\le2\rho\,\kappa^{j}$ uniformly in $p$, so they converge locally uniformly
to $x(\cdot)$, which is therefore holomorphic. Each iterate lies in the closed set
$X_{\mathbb R}$ for $p\in D_{\mathbb R}$, giving the real-subspace claim. The bound follows from $x(p)=T(x(p)),$ since $x(p)$ is a fixed point. By the triangle inequality 
\begin{align*}
\|x(p)-x_0\|
&\le \|x(p)-T(x_0)\|+\|T(x_0)-x_0\| \\
&= \|T(x(p))-T(x_0)\|+\|T(x_0)-x_0\| \\
&\le \kappa\|x(p)-x_0\|+\|T(x_0)-x_0\|
\end{align*}
from the contraction hypothesis. Hence
\[
(1-\kappa)\|x(p)-x_0\|\le\|T(x_0)-x_0\|
\]
which proves the bound.
\end{proof}

\begin{lemma}[Exponential decay of Fourier coefficients on a strip]
\label{lem:strip-decay}
If $\Lambda$ is holomorphic on $S_\delta$ with $\sup_{S_\delta}|\Lambda|\le M$ and
$\Lambda(\theta)=\sum_n c_ne^{in\theta}$ along $\mathbb R$, then $|c_n|\le Me^{-\delta'|n|}$ for
every $\delta'<\delta$ and all $n$.
\end{lemma}

\begin{proof}
By $2\pi$-periodicity of $\Lambda(\theta)e^{-in\theta}$, Cauchy's theorem lets us shift the
defining contour $c_n=\frac1{2\pi}\int_0^{2\pi}\Lambda(\theta)e^{-in\theta}\,d\theta$ to
$\operatorname{Im}\theta=\eta_0$ (the vertical edges cancel), giving
$c_n=e^{n\eta_0}\frac1{2\pi}\int_0^{2\pi}\Lambda(\theta+i\eta_0)e^{-in\theta}\,d\theta$, with
$|c_n|\le Me^{n\eta_0}$. Taking $\eta_0\to\mp\delta'$ with sign opposite to $n$ (Figure 2 shows the
case $n<0$) gives the claim.
\end{proof}

\begin{figure}
\begin{center}
\begin{tikzpicture}[scale=1.1]
  \fill[blue!8] (-0.5,-1.8) rectangle (7.0,1.8);
  \draw[dashed] (-0.5,1.5) -- (7.0,1.5);
  \draw[dashed] (-0.5,-1.5) -- (7.0,-1.5);

  \draw (-0.5,1.5) -- (-0.6,1.5) node[left] {\footnotesize $\delta$};
  \draw (-0.5,-1.5) -- (-0.6,-1.5) node[left] {\footnotesize $-\delta$};
  \draw[red] (-0.5,1.05) -- (-0.6,1.05) node[left] {\footnotesize $\delta'$};
  \draw[red] (-0.5,-1.05) -- (-0.6,-1.05) node[left] {\footnotesize $-\delta'$};
  \draw[red] (-0.5,0.7) -- (-0.6,0.7) node[left,xshift=-2pt] {\footnotesize $\eta_0$};

  \draw[thick,-latex] (0,0) -- (2*3.1416,0);

  \draw[thick,red,-latex] (0,0.7) -- (2*3.1416,0.7);

  \draw[thick,gray] (0,0) -- (0,0.7);
  \draw[thick,gray] (2*3.1416,0) -- (2*3.1416,0.7);

  \draw[red,dotted,-latex] (3.6,0.7) -- (3.6,1.0);

  \draw[->] (-0.5,-2.1) -- (-0.5,2.1) node[above] {\footnotesize $\operatorname{Im}\theta$};
  \node at (0,-0.3) {\footnotesize $0$};
  \node at (2*3.1416,-0.3) {\footnotesize $2\pi$};

  \node[align=center] at (2*3.1416/2,-2.0)
    {\footnotesize $\eta_0\to\mp\delta'$ (sign opposite $n$)};
\end{tikzpicture}
\end{center}
\caption{The contour $\operatorname{Im}\theta=0$ (black) defining $c_n$ is shifted to
$\operatorname{Im}\theta=\eta_0$ (red); the vertical edges at $\theta=0,2\pi$ cancel by
$2\pi$-periodicity of $\Lambda(\theta)e^{-in\theta}$, leaving
$c_n=e^{n\eta_0}\cdot\frac1{2\pi}\int_0^{2\pi}\Lambda(\theta+i\eta_0)e^{-in\theta}\,d\theta$.
Letting $\eta_0\to\delta'$ (shown) gives the bound for $n<0$; the mirrored shift
$\eta_0\to-\delta'$ handles $n>0$.}
\end{figure}
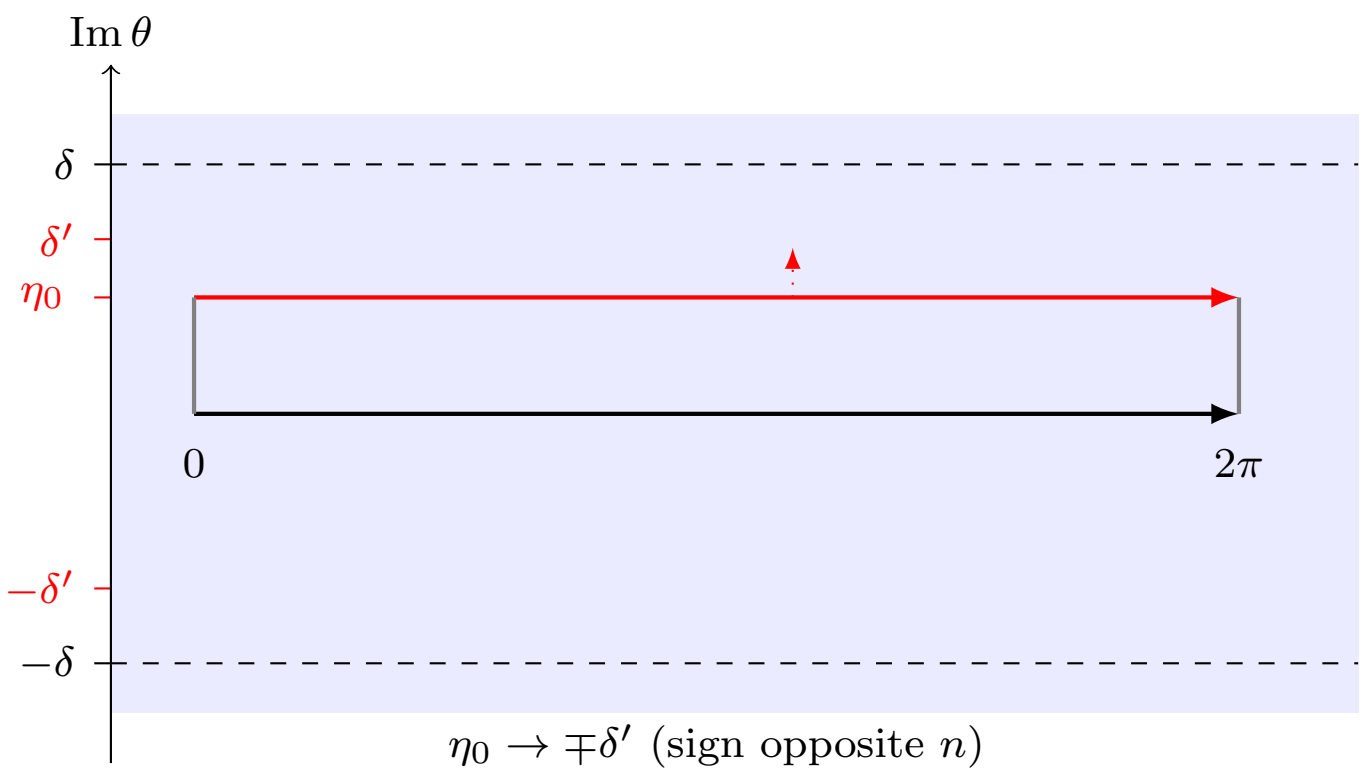

\begin{lemma}[Boundary correspondence for exterior domains]
\label{lem:boundary-converse}
Let $\Gamma\subset\mathbb C$ be a Jordan curve with bounded and unbounded complementary domains
$\operatorname{int}\Gamma$, $\operatorname{ext}\Gamma$, and let
$U:=\{|w|>1\}\cup\{\infty\}$. Suppose $\Phi:\overline U\to\widehat{\mathbb C}$ is holomorphic
and nonconstant on $\{|w|>1\}$ with a simple pole at $\infty$ and $\Phi(w)=cw+O(1)$
($c\ne0$), continuous on $\overline U$, and restricts to a homeomorphism $\{|w|=1\}\to\Gamma$
traversing $\Gamma$ in its positive orientation as $\sigma$ increases in $w=e^{i\sigma}$. Then
$\Phi$ is a conformal bijection of $\{|w|>1\}$ onto $\operatorname{ext}\Gamma$, extending to a
bijection $\overline U\to\overline{\operatorname{ext}\Gamma}$ with $\infty\mapsto\infty$.
\end{lemma}

\begin{proof}
The closed curve $\sigma\mapsto\Phi(e^{i\sigma})$ has winding number
$n(a)$ constant on each component of $\mathbb C\setminus\Gamma$, vanishing on
$\operatorname{ext}\Gamma$; since it traverses $\Gamma$ once positively, $n(a)=1$ on
$\operatorname{int}\Gamma$.

Fix $a\notin\Gamma$. Then $\Phi-a$ is nonvanishing on $|w|=1$ and for
large $|w|$, so for all but finitely many $R>1$ the argument principle on $1<|w|<R$ gives
$N(a,R):=\#\{w:1<|w|<R,\ \Phi(w)=a\}=\frac1{2\pi i}\oint_{|w|=R}\frac{\Phi'}{\Phi-a}
-\frac1{2\pi i}\oint_{|w|=1}\frac{\Phi'}{\Phi-a}$. The second integral is $n(a)$. In the first,
$\Phi(w)/w\to c\ne0$, so for $R$ large the curve $w\mapsto\Phi(w)$ on $|w|=R$ is homotopic in
$\mathbb C\setminus\{a\}$ to $w\mapsto cw$, of winding number $1$; hence $N(a,R)=1-n(a)$ for all
large $R$. As $N(a,R)$ is a nonnegative integer nondecreasing in $R$ and converging to
$1-n(a)\in\{0,1\}$, it satisfies $N(a,R)\le1-n(a)$ for every $R$. Thus $\Phi$ omits every
$a\in\operatorname{int}\Gamma$ on $\{|w|>1\}$, so
\[
\Phi(\{|w|>1\})\cap\operatorname{int}\Gamma=\varnothing,
\tag{5}
\]
and attains every $a\in\operatorname{ext}\Gamma$ exactly once there.

If $\Phi(w_0)=p\in\Gamma$ with $|w_0|>1$, then by the open mapping
theorem\cite{ahlfors1979} $\Phi$ maps a neighborhood of $w_0$ onto a neighborhood of $p$; since
$\Gamma=\partial(\operatorname{int}\Gamma)$, that neighborhood meets $\operatorname{int}\Gamma$,
contradicting (5). Hence $\Phi(\{|w|>1\})\subset\operatorname{ext}\Gamma$, and with the
exact-once count $\Phi$ is a holomorphic bijection, hence conformal, onto
$\operatorname{ext}\Gamma$; with $\Phi(\infty)=\infty$ and the given boundary homeomorphism it
extends to $\overline U\to\overline{\operatorname{ext}\Gamma}$.
\end{proof}

\begin{lemma}[Holomorphy of parameter integrals]
\label{lem:param-hol}
Let $D\subset\mathbb C$ be open, $[a,b]\subset\mathbb R$ compact, and
$F:[a,b]\times D\to\mathbb C$ jointly continuous with $F(\theta,\cdot)$ holomorphic on $D$ for
each fixed $\theta\in[a,b]$. Then $a\mapsto\int_a^bF(\theta,a)\,d\theta$ is
holomorphic on $D$.
\end{lemma}

\begin{proof}
The map is continuous on $D$ by joint continuity of $F$ on the compact set $[a,b]$. For any
triangle $\triangle\subset D$, Fubini's theorem lets us interchange the order of integration in
$\oint_{\partial\triangle}\int_a^bF(\theta,a)\,d\theta\,da$, and Cauchy's
theorem gives $\oint_{\partial\triangle}F(\theta,a)\,da=0$ for each fixed
$\theta$ (as $a\mapsto F(\theta,a)$ is holomorphic), so the double integral
vanishes for every $\triangle\subset D$. Morera's theorem gives the claim.
\end{proof}

\subsection{Sharp starlikeness threshold}

\begin{proposition}
\label{prop:sharp}
For real $a$, $K_a$ is starlike about $0$ if and only if $|a|<\alpha/k$; at
$|a|=\alpha/k$ the curve has a cusp.
\end{proposition}

\begin{proof}
With $N_a:=\operatorname{Im}(z_a'\overline{z_a})
=\alpha^2+ka^2+(k+1)\alpha\cos(mt)$ one has
$\theta_a'=N_a/|z_a|^2$ where $z_a\ne0$. The quadratic
$q(x)=kx^2-(k+1)\alpha x+\alpha^2$ has roots $\alpha/k,\alpha$ and $q(0)>0$. If
$|a|<\alpha/k$ then $N_a\ge q(|a|)>0$ and
$|z_a|\ge\alpha-|a|>0$, so $\theta_a'>0$; since
$|a e^{ikt}|<|\alpha e^{it}|$, Rouché's theorem \cite{ahlfors1979} gives winding number $1$, so
$\theta_a$ is a bijection of $\mathbb T$ and $K_a$ is starlike. If
$\alpha/k\le|a|\le\alpha$ then $q(|a|)\le0$, and at $t$ with
$\cos(mt)=-\operatorname{sgn}a$ we get $N_a\le0$, so $\theta_a$ is
not increasing; if $|a|>\alpha$ the winding number is $k\ge2$, so some ray meets
$K_a$ twice. Finally $z_a'=i\alpha e^{it}(1+\tfrac{ka}\alpha e^{imt})$
vanishes exactly when $k|a|=\alpha$ and $e^{imt}=-\operatorname{sgn}a$,
producing the cusp.
\end{proof}

Since $a_*<\alpha/k$, every real $a\in D_*$ yields a starlike analytic
Jordan curve.

\subsection{Step A: the polar log-radius}

\begin{theorem}
\label{thm:stepA}
For $a\in D_*$:
\begin{enumerate}
\item[\textup{(i)}]  For real $a$, $K_a$ is an analytic Jordan curve,
starlike about $0$, positively oriented, with $0\in\operatorname{int}P_a$; and
$\theta_a|_{\mathbb R}$ is an increasing bijection of $\mathbb T$.
\item[\textup{(ii)}] $\theta_a$ is injective on $S_\delta$ with
$\theta_a(S_\delta)\supseteq S_{3\delta/4}$, and
$s(\theta,a):=\theta_a^{-1}(\theta)$ is holomorphic on
$S_{3\delta/4}\times D_*$.
\item[\textup{(iii)}] $\Lambda_a$ is holomorphic and $2\pi$-periodic on
$S_{\delta/2}$, real on $\mathbb R$ for real $a$, with
\[
\sup_{S_{\delta/2}}|\Lambda_a-\log\alpha|\le\frac{10|a|}{\alpha}.
\tag{6}
\]
\item[\textup{(iv)}] $a\mapsto\Xi_a$ is holomorphic on $D_*$ into
$\mathcal A_{\delta/4}$, with $\|\Xi_a\|_{\delta/8}\le\frac{10|a|}{\alpha}\Sigma_0$.
\end{enumerate}
In particular $\Lambda_0\equiv\log\alpha$ and $\theta_0=\mathrm{id}$.
\end{theorem}

\begin{proof}
\emph{(i)} is Proposition~\ref{prop:sharp} with $a_*<\alpha/k$: regularity is
$z_a'\ne0$ (as $k|a|<\alpha$), positive orientation is
$\theta_a'>0$, and $0\in\operatorname{int}P_a$ since every ray meets
$K_a$ once.

\emph{(ii)} Differentiate (1) directly. With
$h_a=\alpha+a e^{imt}$, $h_a^\#=\alpha+e^{-imt}$,
\[
\theta_a'
=1+\frac m2\Bigl[\frac {a e^{imt}}{h_a}+\frac {a e^{-imt}}{h_a^{\#}}\Bigr]
=1+\frac m2\cdot\frac{\alpha(a e^{imt}+a e^{-imt})+2a^2}{h_a h_a^{\#}}
=\frac{h_a h_a^\#+\tfrac m2\alpha(a e^{imt}+{a e^{-imt}})+ma^2}{h_a h_a^\#}.
\]
Using $h_a h_a^\#=\alpha^2+\alpha({a e^{imt}}+{a e^{-imt}})+a^2$, the numerator is
$\alpha^2+ka^2+\tfrac{(k+1)\alpha}2({a e^{imt}}+{a e^{-imt}})$ (recall $m=k-1$), which agrees with
$\operatorname{Im}(z_a'\overline{z_a})$ on $\mathbb R$. By (4) and
$|a|\le\alpha/(8(k+1))$, a short computation gives
\[
|\text{numerator}-\alpha^2|\le\tfrac{17}{64}\alpha^2,
\qquad
|h_a h^{\#}_a-\alpha^2|\le\tfrac{25}{144}\alpha^2,
\tag{7}
\]
so $|h_a h^{\#}_a|\ge\tfrac{119}{144}\alpha^2$ and
\[
|\theta_a'-1|
=\frac{|\text{numerator}-h_a h^{\#}_a|}{|h_a h^{\#}_a|}
\le\frac{\tfrac{17}{64}+\tfrac{25}{144}}{\tfrac{119}{144}}<\tfrac35
\qquad\text{on }S_\delta .
\tag{8}
\]

Integrating $\theta_a'-1$ along segments in the convex set $S_\delta$,
(8) gives $|\theta_a(t_1)-\theta_a(t_2)|\ge\tfrac25|t_1-t_2|$, hence
injectivity. The constants cancel in (1), so
$\theta_a(t)-t=\operatorname{Im}[\operatorname{Log}(1+{a e^{imt}}/\alpha)-\operatorname{Log}(1+{a e^{-imt}}/\alpha)]$;
since $|\operatorname{Log}(1+w)|\le2|w|$ for $|w|\le\tfrac12$, and $|{a e^{imt}}/\alpha|,|{a e^{-imt}}/\alpha|\le
2|a|/\alpha\le\tfrac12$,
\[
|\theta_a(t)-t|\le\frac{4|a|}\alpha\le\frac\delta8
\qquad\text{on }S_\delta,
\tag{9}
\]
using $a_*\le\alpha\delta/32$. Fix $\theta\in S_{3\delta/4}$ and set
$T(\xi):=-(\theta_a(\theta+\xi)-(\theta+\xi))$ for $|\xi|\le\delta/8$; then
$\theta+\xi\in S_{7\delta/8}\subset S_\delta$, and so $T$ maps $\overline{D(0,\delta/8)}$ into itself
by (9). $|\partial_\xi T|\le\tfrac35$ by (8), which implies that $T$ is a contraction. Thus 
$T$ is jointly holomorphic in $(\xi,\theta,a)$. 
Lemma~\ref{lem:fixed-point} gives a holomorphic fixed point
$\xi(\theta,a)$, and $s:=\theta+\xi$ satisfies $\theta_a(s)=\theta$;
injectivity makes $s=\theta_a^{-1}$, and uniqueness gives
$s(\theta+2\pi,a)=s(\theta,a)+2\pi$.

\emph{(iii)} By (7), $h_a h^{\#}_a\in D(\alpha^2,\tfrac{25}{144}\alpha^2)$,
a disc in the open right half-plane; so the principal root 
$\rho_a=\sqrt{h_a h^{\#}_a}$ is holomorphic,
nonvanishing, with $\operatorname{Re}\rho_a>0$ and $|\rho_a+\alpha|>\alpha$.
Factoring $\rho_a^2-\alpha^2=(\rho_a-\alpha)(\rho_a+\alpha)$,
\[
|\rho_a-\alpha|
=\frac{|h_a h^{\#}_a-\alpha^2|}{|\rho_a+\alpha|}
\le\frac{4\alpha|a|+|a|^2}{\alpha}\le5|a| .
\]
As $5a_*/\alpha<\tfrac12$, $\gamma$, the segment between $\alpha$ and $\rho_a$, lies in
$D(\alpha,\alpha/2)$. For $z\in D(\alpha,\alpha/2)$, $|1/z|\le2/\alpha$, so
\[
|\log\rho_a-\log\alpha|
=|\int_\gamma\frac{1}{z}dz|
\le\tfrac2\alpha|\rho_a-\alpha|
\le\frac{10|a|}{\alpha}
\].
Now $\log\rho_a$ is holomorphic on $S_{7\delta/8}$ and $s=\theta_a^{-1}$
maps $S_{3\delta/4}\times D_*$ holomorphically into $S_{7\delta/8}$ by (ii), so
$\Lambda_a=\log\rho_a\circ s$ is holomorphic on $S_{\delta/2}\times D_*$
(as $\delta/2<3\delta/4$); reparametrizing by $s$ leaves the range of values unchanged, so (6)
follows. $\log\rho_a$ is $2\pi$-periodic and
$s(\theta+2\pi,a)=s(\theta,a)+2\pi$, so $\Lambda_a(\theta+2\pi)=\Lambda_a(\theta)$.
For real $a$, $\rho_a,\theta_a,s$ are real on $\mathbb R$, hence so
is $\Lambda_a$.

\emph{(iv)} By (ii)--(iii), $\Lambda_a(\theta)e^{-in\theta}$ is jointly continuous
on $[0,2\pi]\times D_*$ and holomorphic in $a$ for fixed $\theta$, so
Lemma~\ref{lem:param-hol} gives that $c_n(a):=\widehat{\Xi_a}(n)$ is
holomorphic in $a$. By (6) and Lemma~\ref{lem:strip-decay},
\[
|c_n(a)|\le\frac{10|a|}{\alpha}\,e^{-\frac7{16}\delta|n|}\qquad(n\in\mathbb Z,\ a\in D_*).
\tag{10}
\]
Lemma~\ref{lem:fourier-hol} with $\nu=\tfrac7{16}\delta$ and $\mu=\tfrac14\delta$ gives
holomorphy into $\mathcal A_{\delta/4}$ by (10); with $\mu=\delta/8$ it gives
$\|\Xi_a\|_{\delta/8}\le\frac{10|a|}{\alpha}\Sigma(\tfrac5{16}\delta)\le\frac{10|a|}{\alpha}\Sigma_0$.
At $a=0$, $h_0=h_0^{\#}\equiv\alpha$, so $\rho_0\equiv\alpha$, $\theta_0=\mathrm{id}$,
$\Lambda_0\equiv\log\alpha$.
\end{proof}

\subsection{Step B: the composition operator}

For $\Xi$ holomorphic and bounded on $S_{\delta/2}$ and $u\in\mathcal A$ with
$\|u\|_{\delta/8}\le\delta/16$, define $\mathcal C(\Xi,u)(\sigma):=\Xi(\sigma+u(\sigma))$. This is
well defined: by (2), $\|u\|_\infty\le\delta/16<\delta/2$, so $\sigma+u(\sigma)\in
S_{\delta/2}$ for real $\sigma$. (We allow $u$ complex-valued, as occurs for complex
$a$.)

\begin{theorem}
\label{thm:stepB}
Suppose $|\hat\Xi(n)|\le We^{-\frac7{16}\delta|n|}$ for all $n$, where $W$ is a positive constant, and
$\|u\|_{\delta/8},\|v\|_{\delta/8}\le{\delta/16}$. Define $P_j(\Xi):=\frac{i^j}{j!}\sum_n\hat\Xi(n)n^je^{in\sigma}.$ Then 
\item[\textup{(i)}]$\mathcal C(\Xi,u)=\sum_{j\ge0}P_j(\Xi)u^j$ in
$\mathcal A.$
\item[\textup{(ii)}] $\|P_j(\Xi)\|_{\delta/8}\le W\Sigma_0(\delta/8)^{-j} \qquad (11)$

\item[\textup{(iii)}] $\|\mathcal C(\Xi,u)\|_{\delta/8}\le2W\Sigma_0$

\item[\textup{(iv)}] $\|\mathcal C(\Xi,u)-\mathcal C(\Xi,v)\|_{\delta/8}\le\frac{32W\Sigma_0}{\delta}\|u-v\|_{\delta/8}\qquad(12)$

\item[\textup{(v)}] If $a\mapsto\Xi_a$ has Fourier coefficients holomorphic in $a$
with $|\hat\Xi_a(n)|\le We^{-\frac7{16}\delta|n|}$ uniformly on $D_*$, then
$(u,a)\mapsto\mathcal C(\Xi_a,u)$ is jointly holomorphic on
$\{\|u\|_{\delta/8}\le\delta/16\}\times D_*$; and if $\Xi$ is real-valued, so is
$\mathcal C(\Xi,u)$.
\end{theorem}

\begin{proof}
\emph{(i)} Substituting $\sigma\mapsto\sigma+u(\sigma)$ into
$\Xi=\sum_n\hat\Xi(n)e^{in\sigma}$ and expanding $e^{inu}=\sum_j\frac{(in)^j}{j!}u^j$ gives the
stated double series. It converges absolutely in $\mathcal A$: summing the $j$-series first
(via (2)) contributes a factor $e^{|n|\|u\|_{\delta/8}}$, and since $\|u\|_{\delta/8}\le\delta/16$ 
the remaining exponent is $-\delta/4<0$, so
$\sum_n|\hat\Xi(n)|e^{({\delta/8}+\|u\|_{\delta/8})|n|}\le W\sum_n e^{-\frac14\delta|n|}<\infty$.
Reordering by Fubini's theorem and grouping by powers of $u$ yields $\mathcal C(\Xi,u)=\sum_{j\ge0}P_j(\Xi)u^j$.

\emph{(ii)} From $x^j\le j!\frac{\delta}{8}^{-j}e^{\delta x/8}$ (a single term of
$e^{\delta x/8}$), applied at $x=|n|$,
\[
\|P_j(\Xi)\|_{\delta/8}\le\frac1{j!}\sum_n|\hat\Xi(n)||n|^je^{\delta|n|/8}
\le({\delta/8})^{-j}\sum_n|\hat\Xi(n)|e^{\delta|n|/4}
\le W({\delta/8})^{-j}\sum_ne^{-\frac3{16}\delta|n|}=W\Sigma_0({\delta/8})^{-j},
\]
which proves (11).

\emph{(iii)} With $\|u\|_{\delta/8}\le\delta/16$, (11) gives
$\|\mathcal C(\Xi,u)\|_{\delta/8}\le W\Sigma_0\sum_j2^{-j}=2W\Sigma_0$. 

\emph{(iv)} For the Lipschitz bound,
$u^j-v^j=(u-v)\sum_{r=0}^{j-1}u^rv^{j-1-r}$ gives
$\|u^j-v^j\|_{\delta/8}\le j(\delta/16)^{j-1}\|u-v\|_{\delta/8}$, so
$\sum_{j\ge1}\|P_j\|_{\delta/8}\|u^j-v^j\|_{\delta/8}\le W\Sigma_0\|u-v\|_{\delta/8}\sum_{j\ge1}j({\delta/8})^{-j}(\delta/16)^{j-1} = \frac{16W\Sigma_0}{\delta}\|u-v\|_{\delta/8}\sum_{j\ge1}j2^{-j} =\frac{32W\Sigma_0}{\delta}\|u-v\|_{\delta/8}$
(using $\sum_{j\ge1} j2^{-j}=2$).

\emph{(v)} The $n$-th coefficient of $P_j(\Xi_a)$ is
$\frac{i^j}{j!}\hat\Xi_a(n)n^j$, holomorphic in $a$ and bounded by
$W({\delta/8})^{-j}e^{-\frac5{16}\delta|n|}$ (absorbing $|n|^j$ as above), so
Lemma~\ref{lem:fourier-hol} makes each $a\mapsto P_j(\Xi_a)$ holomorphic
into $\mathcal A$; by (11) the series $\sum_jP_j(\Xi_a)u^j$ converges uniformly on
$\{\|u\|_{\delta/8}\le\delta/16\}\times D_*$ (dominated by $\sum_j W\Sigma_0 2^{-j}$), so its sum
is jointly holomorphic. Finally, if $\Xi$ is real, then $\Xi(\sigma+u(\sigma))$ is real for
each real $\sigma$.
\end{proof}

\subsection{Step C: the boundary correspondence}

\begin{theorem}
\label{thm:stepC}
\item[\textup{(i)}] Define $T_a(u):=-\mathcal H[\mathcal C(\Xi_a,u)]$ and $\overline B:=\{u\in\mathcal A:\|u\|_{\delta/8}\le\delta/16\}$. For
every $a\in D_*$, $T_a$ is a $\tfrac12$-contraction of $\overline B$ into
itself, with a unique fixed point $u_a\in\overline B$. Furthermore, this unique fixed point $u_a$ satisfies 

$u_a+\mathcal H[\mathcal C(\Lambda_a,u_a)]=0, \qquad (13) \qquad \|u_a\|_{\delta/8}\le\frac{20|a|}{\alpha}\Sigma_0\le\tfrac\delta{32}. \qquad (14)$
\item[\textup{(ii)}] $a\mapsto u_a$ is holomorphic on $D_*$ with $u_0=0$,
and $u_a$ is real-valued for real $a$.
\item[\textup{(iii)}] For real $a\in D_*$, $\tau_a:=\sigma+u_a(\sigma)$
is an increasing real-analytic homeomorphism of $\mathbb T$.
\item[\textup{(iv)}] $\tau_a$ is the boundary correspondence of the normalized
exterior conformal map $\Psi_a:\{|w|>1\}\to\Omega_a$; that is,
$\Psi_a(e^{i\sigma})=z_a\bigl(\theta_a^{-1}(\tau_a(\sigma))\bigr)$
for $\sigma\in\mathbb T$. Moreover,
$
\log\operatorname{cap}(K_a)=
\widehat{[\mathcal C(\Lambda_a,u_a)]}(0). \qquad (15)
$
\end{theorem}

\begin{proof}
\emph{(i)} Since
$\mathcal C(\Lambda_a,u)=\mathcal C(\Xi_a,u)+\log\alpha$ and $\mathcal H$
annihilates constants, the fixed-point equation $T_a(u)=u$ is exactly
(13). By (10), $\Xi_a$ meets the hypothesis of
Theorem~\ref{thm:stepB} with $W=\frac{10|a|}{\alpha}$; applying $\mathcal H$ (of norm $\le1$) to
(12) gives, for $u,v\in\overline B$,
\[
\|T_a(u)\|_{\delta/8}\le\frac{20|a|}{\alpha}\Sigma_0,
\qquad
\|T_a(u)-T_a(v)\|_{\delta/8}\le\frac{320|a|\Sigma_0}{\alpha\delta}\|u-v\|_{\delta/8}.
\]
On $D_*$, $\frac{10|a|}{\alpha}\Sigma_0\le\delta/64$, so $\|T_a(u)\|_{\delta/8}\le\delta/32\le\delta/16.$ 
This means $T_a(u)\in\overline B$ whenever $u\in\overline B$ with the Lipschitz constant $\le\tfrac12$. Thus Lemma~\ref{lem:fixed-point} applies
with $\kappa=\tfrac12$, $x_0=0$, $\rho=\delta/16$, giving $u_a$ with
$\|u_a\|_{\delta/8}\le2\|T_a(0)\|_{\delta/8}$; and
$T_a(0)=-\mathcal H[\Xi_a]$ gives
$\|u_a\|_{\delta/8}\le\frac{20|a|}{\alpha}\Sigma_0\le\delta/32$, which is (14).

\emph{(ii)} By Theorem~\ref{thm:stepB},
$(u,a)\mapsto\mathcal C(\Xi_a,u)$ is jointly holomorphic on
$\overline B\times D_*$; composing with the bounded operator $-\mathcal H$, so is
$T_a(u)$, and Lemma~\ref{lem:fixed-point} makes $a\mapsto u_a$
holomorphic on $D_*$. For the real subspace $\mathcal A_{\mathbb R}$ of real-valued functions:
for real $a$, $\Xi_a$ is real on $\mathbb R$, so
$\mathcal C(\Xi_a,u)\in\mathcal A_{\mathbb R}$ whenever $u\in\mathcal A_{\mathbb R}$;
and $\mathcal H$ preserves $\mathcal A_{\mathbb R}$. Hence
$T_a(\overline B\cap\mathcal A_{\mathbb R})\subset\mathcal A_{\mathbb R}$, and
Lemma~\ref{lem:fixed-point} gives $u_a$ real for real $a$. Since
$\Xi_0\equiv0$, $u_0=0$.

\emph{(iii)} By (3) and (14),
$\|u_a'\|_\infty\le\tfrac8{e\delta}\cdot\tfrac\delta{32}=\tfrac1{4e}<1$, so for real
$a$, $\tau_a'=1+u_a'\ge1-\tfrac1{4e}>0$. Thus
$\tau_a=\mathrm{id}+u_a$ is a strictly increasing, $2\pi$-periodic-mod-$2\pi$,
real-analytic bijection of $\mathbb T$, i.e.\ an increasing homeomorphism.

\emph{(iv)} Fix real $a\in D_*$ and set
$\Theta:=\mathcal C(\Lambda_a,u_a)=\Lambda_a\circ\tau_a\in\mathcal A$,
real-valued, with $\Upsilon:=\widehat\Theta(0)\in\mathbb R$. Since $\sum_n|\widehat\Theta(n)|<\infty$,
\[
G(w):=\sum_{n\ge1}2\widehat\Theta(-n)w^{-n},
\qquad
\Phi(w):=e^{\Upsilon}we^{G(w)}
\]
define $G$ holomorphic on $\{|w|>1\}$, continuous on $\{|w|\ge1\}$, with $G(\infty)=0$; the
absence of a constant term in $G$ is what makes the boundary computation below exact.

As $\Theta$ is real, $\widehat\Theta(n)=\overline{\widehat\Theta(-n)}$.
Comparing Fourier series term by term on $w=e^{i\sigma}$,
\[
\operatorname{Re}G(e^{i\sigma})
=\operatorname{Re}\sum_{n\ge1}2\widehat\Theta(-n)e^{-in\sigma}
=\sum_{n\ge1}\bigl(\widehat\Theta(n)e^{in\sigma}+\widehat\Theta(-n)e^{-in\sigma}\bigr)
=\sum_{n\ne0}\widehat\Theta(n)e^{in\sigma}
=\Theta(\sigma)-\Upsilon,
\]
\[
\operatorname{Im}G(e^{i\sigma})
=\operatorname{Im}\sum_{n\ge1}2\widehat\Theta(-n)e^{-in\sigma}
=i\sum_{n\ne0}\operatorname{sgn}(n)\,\widehat\Theta(n)e^{in\sigma}
=-\mathcal H[\Theta](\sigma)
=u_a(\sigma),
\]

the last equality by (13). (Both series converge absolutely, so these are genuine identities on
$\mathbb T$.) Hence on $|w|=1$,
\[
\Phi(e^{i\sigma})=e^{i\sigma}e^{\Upsilon+\operatorname{Re}G+i\operatorname{Im}G}
=e^{\Theta(\sigma)}e^{i(\sigma+u_a(\sigma))}
=e^{\Lambda_a(\tau_a(\sigma))}e^{i\tau_a(\sigma)}
=z_a\bigl(\theta_a^{-1}(\tau_a(\sigma))\bigr),
\tag{16}
\]
using $z_a=\rho_a e^{i\theta_a}$ and
$\Lambda_a=\log\rho_a\circ\theta_a^{-1}$. Thus every value lies on
$K_a$. As $w\to\infty$, $G(w)\to0$, so $\lim_{w\to\infty}\frac{\Phi(w)}{w}=e^{\Upsilon}e^{G(w)}=e^{\Upsilon}$. Hence $\Phi$ extends holomorphically to $\infty$ with a simple pole, 
$\Phi(\infty)=\infty$, and leading coefficient $e^{\Upsilon}>0$.

By Theorem~\ref{thm:stepA}(i), $\theta_a^{-1}$ is an increasing homeomorphism of $\mathbb T$;
so is $\tau_a$ by (iii), hence so is their composition. By (16), $\sigma\mapsto\Phi(e^{i\sigma})$
is therefore an increasing homeomorphism of $\mathbb T$ onto the Jordan curve $K_a$.
We check the remaining hypotheses of Lemma~\ref{lem:boundary-converse}. First, $\Phi$ is continuous on
$\{|w|\ge1\}$ because $\sum_n|\widehat\Theta(n)|<\infty$, and $\Phi(w)=e^{\Upsilon}w+O(1)$ with
$e^{\Upsilon}\neq0$ because $e^{G(w)}=1+O(1/w)$. Second, by (16),
$\Phi(e^{i\sigma})=e^{\Theta(\sigma)}e^{i\tau_a(\sigma)}$ never vanishes and
$\tau_a(\sigma+2\pi)=\tau_a(\sigma)+2\pi$, so the winding number of $\Phi(e^{i\sigma})$ about $0$ is $+1$.
Since $K_a$ is a Jordan curve, its winding number is constant on $\operatorname{int}K_a$ and zero on
$\operatorname{ext}K_a$; as it is nonzero about $0$, we get $0\in\operatorname{int}K_a$, and the
winding number is $+1$ about every point of $\operatorname{int}K_a$.
Taking $\Gamma=K_a$, $\operatorname{int}\Gamma=\operatorname{int}P_a$,
$\operatorname{ext}\Gamma=\Omega_a$, Lemma~\ref{lem:boundary-converse} shows that $\Phi$ is a conformal
bijection of $\{|w|>1\}$ onto $\Omega_a\setminus\{\infty\}$ with $\Phi(w)=e^{\Upsilon}w+O(1)$ and
$e^{\Upsilon}>0$: exactly the normalization of $\Psi_a$. By uniqueness of the normalized exterior
map~\cite[\S1.2]{pommerenke1975}, $\Phi=\Psi_a$ and $e^{\Upsilon}=\operatorname{cap}(K_a)$, so $\tau_a$ is its
boundary correspondence, and $\Upsilon=\widehat\Theta(0)$ gives (15).
\end{proof}

\subsection{Step D: analyticity of the capacity}

\begin{theorem}
\label{thm:stepD}
\[
L(a):=\frac1{2\pi}\int_0^{2\pi}\Lambda_a\bigl(\sigma+u_a(\sigma)\bigr)\,d\sigma
\]
is holomorphic on $D_*$ and equals $\log\operatorname{cap}(K_a)$ for real
$a\in D_*$. Hence $a\mapsto\operatorname{cap}(K_a)=e^{L(a)}$
is real-analytic on $(-a_*,a_*)$, with $\operatorname{cap}(K_0)=\alpha$.
\end{theorem}

\begin{proof}
Fix $\sigma$. By Theorem~\ref{thm:stepC}, $a\mapsto u_a(\sigma)$ is
holomorphic on $D_*$, with $|u_a(\sigma)|\le\delta/32<\delta/2$. Since
$(\theta,a)\mapsto\Lambda_a(\theta)$ is holomorphic on $S_{\delta/2}\times
D_*$ (Theorem~\ref{thm:stepA}), the integrand $\Lambda_a(\sigma+u_a(\sigma))$
is jointly continuous on $[0,2\pi]\times D_*$ and holomorphic in $a$ for fixed
$\sigma$, so Lemma~\ref{lem:param-hol} gives holomorphy of $L$ on $D_*$. For real $a$
the integral is $\widehat\Theta(0)=\log\operatorname{cap}(K_a)$ by (15), so
$\operatorname{cap}(K_a)=e^{L(a)}$ for real $a\in D_*$. $L$ is real there, hence real-analytic, and $\exp$ is entire, so
$a\mapsto\operatorname{cap}(K_a)=e^{L(a)}$ is real-analytic on
$(-a_*,a_*)$.
At $a=0$, $\Lambda_0\equiv\log\alpha$
and $u_0=0$ give $L(0)=\log\alpha$, so $\operatorname{cap}(K_0)=e^{L(0)}=\alpha$.
\end{proof}

\section{Analytic Coefficients of $\gamma(K_a)$}

Throughout this section we specialize the family of Section~1 to the unit circle ($\alpha=1$), $z_a(t):=e^{it}+ae^{ikt}$. When $k=1$ the
perturbation rescales the radius, so $\gamma(K_a)=1+a$; for $k=0$ the perturbation only translates the figure, so $\gamma(K_a)=1$. When $k\le0$ we show below that,
for $|a|$ below an explicit threshold, the perturbation has no effect at all:
$\gamma(K_a)=1$. $k=2$ is of particular importance, since we are able to express $\gamma(K_a)$ as a closed-form expression in terms of $a$. This allows us to solve for the capacity of a cardioid, which is when $a=\pm1/2$, as well as show that all Taylor coefficients in this case are integers. We conjecture that every Taylor coefficient is an integer for all integer values of $k$.

\begin{lemma}[Evenness of the capacity]
\label{lem:even}
Let $\alpha>0$, $k\ge2$, $z_a(t)=\alpha e^{it}+a e^{ikt}$, and
$K_a:=z_a([0,2\pi])$. Then $a\mapsto\operatorname{cap}(K_a)$
is even on any interval of $a$ for which $K_a$ is a Jordan curve.
\end{lemma}

\begin{proof}
Let $c:=\pi/(k-1)$, so $e^{i(k-1)c}=-1$ and hence $e^{ikc}=-e^{ic}$. Then
\[
z_a(t+c)=\alpha e^{ic}e^{it}+a e^{ikc}e^{ikt}
=e^{ic}\bigl(\alpha e^{it}-a e^{ikt}\bigr)=e^{ic}z_{-a}(t).
\]
Since $z_a$ is $2\pi$-periodic, $K_a=\{z_a(t+c):t\in[0,2\pi]\}=e^{ic}K_{-a}$, so
$K_a$ and $K_{-a}$ differ only by rotation. Logarithmic capacity is invariant under rotations\cite{ransford1995}, since these preserve all pairwise distances; hence
$\operatorname{cap}(K_a)=\operatorname{cap}(K_{-a})$.
\end{proof}

\begin{corollary}[Taylor expansion in even powers]
\label{cor:even-powers}
By Theorem 4, $a \mapsto \gamma(K_a)$ is real-analytic at 0. Hence there exists $a_*>0$ such that $\gamma(K_a)=\sum_{m\ge0}\gamma_{2m}a^{2m}$ for $|a|<a_*$, where each
$\gamma_{2m}$ is real, $\gamma_0=1$, and every odd-order coefficient vanishes.
\end{corollary}

\begin{proof}
By Theorem~\ref{thm:stepD}, $a\mapsto\gamma(K_a)$ is real-analytic at $0$, giving a convergent
expansion $\gamma(K_a)=\sum_n\gamma_na^n$ for $|a|<a_*$. By Lemma~\ref{lem:even} (with
$\alpha=1$), $\gamma(K_a)=\gamma(K_{-a})$; substituting $a=-a$ yields a new power series for $\gamma(K_{-a})$. Matching coefficients (unique by uniqueness of power series) gives $\gamma_n=(-1)^n\gamma_n$, so $\gamma_n=0$ for odd $n$.
Finally, $\gamma_0=\gamma(K_0)=\gamma(\overline{\mathbb D})=1$.
\end{proof}

\begin{lemma}[Global formula for the phase derivative]
\label{lem:phase}
Let $\alpha>0$, $k\in\mathbb Z$, $m:=k-1$, and for real $t$ and $|a|<\alpha$ rewrite $z_a(t):=e^{it}+ae^{ikt}$ as
$z_a(t)=\rho_a(t)e^{i\theta_a(t)}$, with $\rho_a=|z_a|$,
$\theta_a=\arg z_a$. Then, uniformly in $t$,
\[
\theta_a'(t)=1+m\sum_{n=1}^\infty(-1)^{n+1}\Bigl(\frac a\alpha\Bigr)^n\cos(nmt).
\]
\end{lemma}

\begin{proof}
Factor $z_a(t)=\alpha e^{it}g(t)$, $g(t):=1+(a/\alpha)e^{imt}$. Since
$|(a/\alpha)e^{imt}|=|a/\alpha|<1$, the triangle inequality gives
$\operatorname{Re}g(t)\ge1-|a/\alpha|>0$ for every $t$, so $\theta_a(t)=t+\arg g(t)$
(principal branch) and, using $\tfrac{d}{dt}\arg g=\operatorname{Im}(g'/g)$ with
$g'=im(a/\alpha)e^{imt}$ and $\operatorname{Im}(iz)=\operatorname{Re}(z)$,
\[
\theta_a'(t)=1+m\operatorname{Re}\!\left(\frac{(a/\alpha)e^{imt}}{1+(a/\alpha)e^{imt}}\right). \tag{17}
\]
Writing $w:=(a/\alpha)e^{imt}$, $|w|=|a/\alpha|<1$, the geometric series
$w/(1+w)=\sum_{n\ge1}(-1)^{n+1}w^n$ converges absolutely, uniformly in $t$ (each term is bounded
by $|a/\alpha|^n$, independent of $t$, so the Weierstrass $M$-test applies); taking
real parts termwise gives
\[
\operatorname{Re}\Bigl(\frac w{1+w}\Bigr)=\sum_{n=1}^\infty(-1)^{n+1}\Bigl(\frac a\alpha\Bigr)^n\cos(nmt),
\]
and substituting yields the above formula for $\theta_a'(t)$.
\end{proof}

\subsection{Coefficients for $k\le1$}

The aim of this section is to prove what was stated in the introduction for section 3, namely that for $|a|$ below an explicit threshold, the $k<1$ perturbation has no effect on analytic capacity.

\begin{lemma}[Phase derivative at negative frequencies]
\label{lem:phase-bound-neg}
Let $k=-n$ with $n\ge1$ and $a\in\mathbb R$, $|a|<1$. Then
\[
\theta_a'(t)=\frac{1-(n-1)a\cos((n+1)t)-na^2}{1+2a\cos((n+1)t)+a^2},
\qquad
\min_{t\in\mathbb R}\theta_a'(t)=\frac{1-n|a|}{1+|a|},
\]
which is strictly positive if and only if $|a|<1/n$.
\end{lemma}

\begin{proof}
In Lemma~\ref{lem:phase} we derived (17) for every integer $k$ with $|a|<1$, giving us
$\theta_a'(t)=1+m\operatorname{Re}\bigl(\tfrac w{1+w}\bigr)=\operatorname{Re}
\bigl(\tfrac{1+kw}{1+w}\bigr)$ for $\alpha=1$, $w:=ae^{imt}$, and $m:=k-1$. Setting $k=-n$,
$w=ae^{-i(n+1)t}$, and multiplying both numerator and denominator by
$\overline{1+w}=1+\overline{w}$, we arrive at
\[
\theta_a'(t)=
\operatorname{Re}\Bigl(\frac{1+kw}{1+w}\cdot\frac{1+\overline{w}}{1+\overline{w}}\Bigr) =
\frac{\operatorname{Re}(1+kw+\overline{w}+kw\overline{w})}{\operatorname{Re}(1+w+\overline{w}+w\overline{w})}
=\frac{1-(n-1)a\cos((n+1)t)-na^2}{1+2a\cos((n+1)t)+a^2},
\]
as stated above.

Write
$x:=\cos((n+1)t)\in[-1,1]$ and $f(x):=\frac{1-(n-1)ax-na^2}{1+2ax+a^2}$,
so that $\theta_a'(t)=f(\cos((n+1)t))$; the denominator satisfies
$1+2ax+a^2\ge1-2|a|+a^2=(1-|a|)^2>0$ for $x\in[-1,1]$, so $f$ is smooth there. By the chain
rule,
\[
\frac{d}{dt}f(\cos((n+1)t))=f'(\cos((n+1)t))\cdot\bigl(-(n+1)\sin((n+1)t)\bigr),
\]
and a direct computation of $f'(x)$ gives
\[
\frac{d}{dt}\theta_a'(t)=\frac{d}{dt}f(\cos((n+1)t))=\frac{a(n+1)^2(1-a^2)\sin((n+1)t)}{(1+2a\cos((n+1)t)+a^2)^2},
\]
which vanishes exactly when $\sin((n+1)t)=0$, i.e.
$x=\cos((n+1)t)=\pm1$.

Since $f'(x)=\dfrac{-(n+1)a(1-a^2)}{(1+2ax+a^2)^2}$ has constant sign $-\operatorname{sgn}(a)$
on $[-1,1]$ (the denominator is a square, hence positive, and $(n+1)>0$, $1-a^2>0$), $f$ is
strictly monotonic there, so its extrema occur at $x=\pm1$, matching the critical points
above; thus $\min_t\theta_a'(t)=\min\{f(1),f(-1)\}$.
Using $1-(n-1)a-na^2=(1-na)(1+a)$, $1+2a+a^2=(1+a)^2$, and the analogous factorizations at
$x=-1$, $f(1)=\frac{1-na}{1+a}, f(-1)=\frac{1+na}{1-a}$.
Therefore
\[
\min_{t\in\mathbb R}\theta_a'(t)=\min\{f(1),f(-1)\}=\frac{1-n|a|}{1+|a|},
\]
positive if and only if $|a|<1/n$.
\end{proof}

\begin{proposition}[Capacity under general negative-frequency perturbations]
\label{prop:neg-freq-general}
Let $a\in\mathbb R$ and $(\beta_k)_{k\ge1}$ satisfy $\sum_{k\ge1}|\beta_k|<\infty$, and suppose
\[
\Gamma_a(t):=e^{it}+a\sum_{k=1}^\infty\beta_ke^{-ikt},\qquad t\in[0,2\pi),
\]
is a positively-oriented parametrization of a Jordan curve $K_a$. Then $\gamma(K_a)=1$. 
\end{proposition}

\begin{proof}
Set $\Phi(z):=z+a\sum_{k\ge1}\beta_kz^{-k}$ for $|z|\ge1$. Since $\sum_k|\beta_k|<\infty$, the
series converges absolutely and uniformly on $|z|\ge1$ (Weierstrass $M$-test),
so $\Phi$ is holomorphic and nonconstant on $\{|z|>1\}$, continuous on $\{|z|\ge1\}$, with a
simple pole at $\infty$. This is the normalization $c=1$, meaning that $\Phi(z)=cz+O(1)$ with $c=1$. On $|z|=1$,
$\Phi(e^{it})=\Gamma_a(t)$, so by hypothesis $\Phi$ restricts to a positively-oriented
homeomorphism $\{|z|=1\}\to K_a$. Lemma~\ref{lem:boundary-converse} therefore identifies $\Phi$
with the normalized exterior conformal map of $K_a$, with leading coefficient $1$; hence
$\gamma(K_a)=1$.
\end{proof}

Taking $\beta_n=1$ and $\beta_j=0$ for $j\ne n$ recovers the single-parameter perturbation $z_a(t)=e^{it}+ae^{-int}$.

\begin{proposition}[The single-parameter case satisfies the hypotheses of Proposition~\ref{prop:neg-freq-general}]
\label{prop:single-parameter}
Let $n\ge1$, $\beta_n=1$, $\beta_j=0$ for $j\ne n$, and $|a|<1/n$. Then
$z_a(t):=e^{it}+ae^{-int}$ parametrizes a Jordan curve $K_a$ in positive orientation, so
Proposition~\ref{prop:neg-freq-general} applies (with $\sum_k|\beta_k|=1<\infty$) and gives
$\gamma(K_a)=1$.
\end{proposition}

\begin{proof}
By Lemma~\ref{lem:phase-bound-neg}, $\theta_a'(t)>0$ for all $t$ under the hypothesis
$|a|<1/n$. Since $|ae^{-int}|=|a|<1=|e^{it}|$, Rouché's theorem\cite{ahlfors1979} gives that
$z_a(t)$ has winding number $1$ about $0$, so $\theta_a(t+2\pi)=\theta_a(t)+2\pi$; together
with $\theta_a'>0$, $t\mapsto\theta_a(t)$ is an increasing bijection of
$\mathbb R/2\pi\mathbb Z$. Since also $|z_a(t)|\ge1-|a|>0$, the map $t\mapsto z_a(t)$ is
injective, so it is a positively-oriented parametrization of a Jordan curve $K_a$. The
hypotheses of Proposition~\ref{prop:neg-freq-general} are met (with $\beta_n=1$, all other
$\beta_j=0$, so $\sum_k|\beta_k|=1<\infty$), giving $\gamma(K_a)=1$.
\end{proof}

For $k=1$, $z_a(t)=(1+a)e^{it}$, so $K_a$ is the circle of radius $|1+a|$. Hence $\gamma(K_a)=1+a$ for $|a|<1$.

\subsection{Closed Form of Capacity when $k=2$}
In this section, we will derive the explicit formula for $\gamma(e^{it}+ae^{2it})$ in terms of $a$, giving us a basis for the exploration of further powers of $k$ and the interesting properties related to this family of sets. 

\begin{lemma}[\cite{ransford1995}]
    Let $E$ be a compact set of $\mathbb{C}$ and let $p(z)=\sum_{k=0}^n a_k z^k$ where $a_n \neq 0$. If 
    \[E_p=\{z \in \mathbb{C}: p(z) \in E\},\]
    then 
    \[\operatorname{Cap(E_p)}=\left(\frac{\operatorname{Cap}(E)}{|a_n|}\right)^{1/n}\]
    where $\operatorname{Cap}(E)$ represents the logarithmic capacity of the set $E$.
    \label{lem: polypre}
\end{lemma}

With this, we can prove the following theorem regarding the analytic capacity of sets bounded by the Jordan curve $e^{it}+ae^{2it}$.

\begin{theorem}
\label{thm: k=2 full function}
    For $a \in (-\frac{1}{2},\frac{1}{2}) \setminus \{0\}$,
    \[\gamma(e^{it}+ae^{2it})=\frac{1-4a^2}{a}\rho(a) \left(\psi(\rho(a)^2)\right)^4=\frac{1-4a^2}{a}\rho(a)\prod_{n=1}^\infty \left(\frac{1-\rho(a)^{4n}}{1-\rho(a)^{4n-2}}\right)^4\]
    where 
    $\rho(a)=\frac{1-\sqrt{1-4a^2}}{2a}$ and $\psi(q)=\sum_{n=1}^{\infty}q^{n(n+1)/2}$ is the Ramanujan Psi function.
\end{theorem}
\begin{proof}
    Note that for $0<|a|<\frac{1}{2}$, the closed Jordan region $K_a $ bounded by the curve $\{e^{it}+ae^{2it}: t\in [0,2\pi]\}$ is path connected. Therefore, we have that $\operatorname{Cap}(K_a)=\gamma(K_a)$. We can also note that the pre-image of the set $K_a$ under the polynomial $p(z)=az^2+z$ is the set 
    \[E_a=\left(\overline{\mathbb{D}}-\frac{1}{a}\right)\cup \overline{\mathbb{D}}\]
    which has the same logarithmic  capacity as 
    \[F_a=\left(\overline{\mathbb{D}}-\frac{1}{2a}\right)\cup \left(\overline{\mathbb{D}}+\frac{1}{2a}\right)\]
    
    Then, Lemma \ref{lem: polypre} gives us that $\gamma(K_a)=\operatorname{Cap(K_a)}=|a|\operatorname{Cap}(E_a)^2=|a|\operatorname{Cap}(F_a)^2$.
    Now, without loss of generality, we can assume that $a>0$. Then from Example 4.6 of  \cite{LiesenSeteNasser2017}, we have that 
    \[\operatorname{Cap}(F_a)=\frac{2K}{\pi }\left(\frac{\sqrt{1-4a^2}}{2a}\right)\sqrt{2L(1+L^2)}\]
    Where $\rho=\frac{1-\sqrt{1-4a^2}}{2a}$, $L=\frac{\theta_2(\rho^4)}{\theta_3(\rho^4)}$, and $K=K(L^2)$ which is the complete elliptic integral of the first kind with elliptic modulus $L^2$. Further, $\theta_2(\cdot)$ and $\theta_3(\cdot)$ are the second and third Jacobi Theta functions respectively.  We can use the fact that $\theta_3^2(\rho^4)=\frac{2K(L^2)}{\pi }$ and the fact that $L=\frac{\theta_2(\rho^4)}{\theta_3(\rho^4)}$ in order to get that 
    \[\operatorname{Cap}(F_a)=\left(\frac{\sqrt{1-4a^2}}{2a}\right)\theta_3^2(\rho^4)\sqrt{\frac{2\theta_2(\rho^4)}{\theta_3(\rho^4)}\left(\frac{\theta_3^2(\rho^4)+\theta_2^2(\rho^4)}{\theta_3^2(\rho^4)}\right)}.\]
    Here we have that $\theta_3(\rho^2)=\theta_3^2(\rho^4)+\theta_2^2(\rho^4)$ and so
    \[\gamma(K_a)=a\left(\frac{1-4a^2}{4a^2}\right)\theta_3^4(\rho^4)\left(\frac{2\theta_2(\rho^4)\theta_3^2(\rho^2)}{\theta_3^3(\rho^4)}\right)=
    \frac{1-4a^2}{4a}(2\theta_2(\rho^4)\theta_3(\rho^4))\theta_3^2(\rho^2).\]
    From Chapter 21 in \cite{WhittakerWatson1927}, we have that $\theta_2^2(\rho^2)=2\theta_2(\rho^4)\theta_3(\rho^4)$, so 
    \[\gamma(K_a)=\frac{1-4a^2}{4a}\left(\theta_2(\rho^2)\theta_3(\rho^2)\right)^2.\]
    We can simplify this further using the product representation of both of these theta functions from the same chapter in \cite{WhittakerWatson1927} in order to get that 
    \[\gamma(K_a)=\frac{1-4a^2}{4a}\left(\rho^{1/2} \prod_{n=1}^\infty  (1-\rho^{4n})^2(1+\rho^{4n})^2(1+\rho^{4n-2})^2\right)^2\]

    From here, we can note that $\prod_{n=1}^\infty(1+\rho^{4n})^2$ is a product of $(1+\rho^k)^2$ where $k$ is a multiple of four and $\prod_{n=1}^\infty(1+\rho^{4n})^2$ is a product of $(1+\rho^k)^2$ where $k$ is a multiple of $2$ but not four, and so $\prod_{n=1}^\infty(1+\rho^{4n})^2(1+\rho^{4n-2})^2=\prod_{n=1}^\infty(1+\rho^{2n})^2$
    where the product goes over all multiples of two. Hence, we get that our product reduces to
    \[=\rho \frac{1-4a^2}{4a}\left( \prod_{n=1}^\infty  (1-\rho^{4n})^2(1+\rho^{2n})^2 \right)^2=\rho \frac{1-4a^2}{4a}\left( \prod_{n=1}^\infty  (1-\rho^{4n})^2\left(\frac{1-\rho^{4n}}{1-\rho^{2n}}\right)^2 \right)^2\]
    

    \[=\frac{1-4a^2}{4a}\rho  \prod_{n=1}^\infty \left(\frac{1-\rho^{4n}}{1-\rho^{4n-2}}\right)^4\]
    as needed. Further, according to part $(ii)$ of entry $21 $ in \cite{berndt1991ramanujan}, we have that this function equals $\frac{1-4a^2}{a}\rho(\psi(\rho^2))^4$.

\end{proof}

We can also solve for the value of a cardioid (the $|a|=\frac{1}{2}$ case) directly. 
\begin{corollary}
    The analytic capacity of a cardioid is $\gamma(e^{it}+\frac{1}{2}e^{2it})=\frac{\pi^2}{8}$.
\end{corollary}
\begin{proof}
    We can continue the same argument as in Theorem \ref{thm: k=2 full function}. Let our region be represented by $K_{1/2}$. In this case, the pre-image of $K_{1/2}$ under $p(z)=\frac{1}{2}z^2+z$ becomes $E_{1/2}=(\overline{\mathbb{D}}-2)\cup \overline{\mathbb{D}}$ which corresponds to two tangent unit discs. So, as before, we have $\gamma(K_{1/2})=\operatorname{Cap}(K_a)=\frac{1}{2}\operatorname{Cap}^2(E_{1/2})$. From \cite{Chen2020TangentDiscs}, we know that $\operatorname{Cap}(E_{1/2})=\frac{\pi}{2}$.
    Combining this, we get that 
    \[\gamma(K_{1/2})=\frac{1}{2}(\frac{\pi}{2})^2=\frac{\pi^2}{8}.\]
\end{proof}

\subsection{Integrality of the Taylor coefficients when $k=2$}

With an explicit expression for $\gamma(e^{it}+ae^{2it})$ in Theorem \ref{thm: k=2 full function}, we can observe a very interesting property regarding its power series expansion with respect to $a$. 

\begin{proposition}
    Let $\gamma(e^{it}+ae^{2it})=\sum_{n=0}^\infty \beta_na^n$ for $0<|a|<\frac{1}{2}$. Then $\beta_n \in \mathbb{Z}$ for all $n \in \mathbb{Z}^{\geq 0}$.
    \label{prop:2int}
\end{proposition}

\begin{proof}
    Let $\mathbb{Z}[[x]]$ represent the ring of formal power series with integer coefficients. Then, we know that $\psi(x)=\sum_{n=0}^\infty x^{n(n+1)/2} \in \mathbb{Z}[[x]]$. We can also show that $\rho(x) \in \mathbb{Z}[[x]]$. To do so, note that $\frac{1-\sqrt{1-4y}}{2y} $ is the generating function for the Catalan numbers, which is an integer sequence. So, $\frac{1-\sqrt{1-4x^2}}{2x^2} \in \mathbb{Z}[[x]]$. Then, lastly, we get that $\rho(x)=x(\frac{1-\sqrt{1-4x^2}}{2x^2})\in \mathbb{Z}[[x]]$. Further, we know that the constant term of $\rho(x)$ is $0$. So, this means that $\rho(x)^2 \in \mathbb{Z}[[x]]$ with zero constant term. So, $\psi(\rho(x)^2)$ and $\psi^4(\rho^2(x))\in \mathbb{Z}[[x]]$. Then we know $1-4x^2$ and $\rho(x)/x$ are in $\mathbb{Z}[[x]]$. So, by closure, we have that $\gamma(e^{it}+xe^{2it})\in \mathbb{Z}[[x]]$.
    
\end{proof}

\begin{conjecture}[Integrality of every Taylor coefficients]
\label{conj:integrality}
For $z_a(t)=e^{it}+ae^{ikt}$ and any integer $k$, every Taylor coefficient $\gamma_j$ of
$\gamma(K_a)$ is an integer.
\end{conjecture}














\bibliographystyle{plain}
\bibliography{references}

\section*{Statement of AI Use}

We used generative AI extensively throughout their research. ChatGPT and Claude were both used to investigate and generate methods of proof. However, wherever AI was used, its results were evaluated based on the comprehension of human authors. We take full responsibility for all content herein.

\end{document}